%Paper: 9205003
%From: wux@milo.math.scarolina.edu (Xian Wu)
%Date: Thu, 7 May 92 15:52:33 EDT

\documentstyle{amsppt}
\magnification=\magstep 1

%convenient abbreviations
\predefine\dotunder{\d}
\redefine\d{\delta}
\predefine\czech{\v}
\redefine\v{\nu}
\define\bc#1{\pmb{\Cal#1}}
\TagsOnRight
\NoBlackBoxes

\topmatter

\title
Chern classes and degenerations of hypersurfaces and their lines
\endtitle

\author Xian Wu
\endauthor

\affil
Department of Mathematics\\
University of South Carolina\\
Columbia, SC 29208
\endaffil

\endtopmatter
\document

\smallpagebreak
\heading 1. Introduction
\endheading
\smallpagebreak

Let $X$ be a generic hypersurface of degree $d$ in $\Bbb P^n$
and $G(2,n+1)$ be the Grassmannian of the lines in $\Bbb P^n$.
The class of the lines contained in $X$
is then equal to the top Chern class of the vector bundle $S^dU^*$,
where $U$ is the universal bundle on $G(2,n+1)$.
Using this one can easily compute that,
for example,
there are $27$ lines on a generic cubic surface
and $2,875$ lines on a generic quintic threefold.
What happens when $X$ degenerates?
For example,
if $X$ becomes reducible,
how can one identify the lines in each component which are the limits
of the lines in a generic $X$ and compute the class of those lines
in the Chow ring of $G(2,n+1)$?
In the case of quintic threefold,
Katz [K] has given a complete description for the geometry of the
limiting lines.
In this paper we will investigate such problems in general.

Let $\{X_s=KL+sD+\cdots \}_{s\in \Delta}$
be a degeneration of generic hypersurfaces of degree $d$
in $\Bbb P^n$,
where $\deg K=k$, $\deg L=l$, and $\deg D=d=k+l$.
We will use the same letter to denote a hypersurface and its
equation.
Define
$$\sigma_K=\{\alpha \in G(2,n+1)| \  \alpha \subset K, \alpha \cap L \subset
D\},$$
$$\sigma_L=\{\alpha \in G(2,n+1)| \  \alpha \subset L, \alpha \cap K \subset
D\}.$$
Katz's results [K] suggest that one should study $\sigma_K$ and $\sigma_L$
in order to understand the set of limiting lines.
By tensoring the section of $S^lU^*$ induced by $L$,
we get a homomorphism $\rho_L$ from
$S^kU^*$ to $S^dU^*$ and
hence have an exact sequence of coherent sheaves
$$S^kU^* @>\rho_L>> S^dU^* @>>> B_K @>>> 0 ,$$
where
$$B_K=S^dU^*/\rho_L(S^kU^*)$$
is the quotient.
Next,
one can formally define the Chern class of $B_K$
by usual multiplicative formula for the quotient bundle;
that is,
one defines inductively
$$\left\{ \eqalign{c_0(B_K) & =1,\cr
c_i(B_K) & =c_i(S^dU^*)-\sum_{j=0}^{i-1} c_j(B_K)c_{i-j}(S^kU^*),
\quad i > 0 .\cr}\right.$$
Similarly, one defines $B_L$ and its class.
Our first result is the following:

\proclaim{Theorem 2.3}
Let
$$c(d,k)=c_{k+1}(S^kU^*)c_{l}(B_K),\quad c(d,l)=c_{l+1}(S^lU^*)c_{k}(B_L).$$
Then,
$c(d,k)$ and $c(d,l)$ are equal to the class of $\sigma_K$ and
the class of $\sigma_L$,
respectively.
\endproclaim

Intuitively,
we may think that the rank of $B_K$ is equal to
$$\text{rank}(S^dU^*)-\text{rank}(S^kU^*)=d+1-(k+1)=l.$$
Since the top Chern class of a vector bundle is equal to the
zero locus of a regular section of the bundle,
one sees that
$$c(d,k)=c_{top}(S^kU^*)c_{``top"}(B_K)$$
does give the kind of the conditions which define $\sigma_K$.
The problem is that $B_K$ is not a quotient bundle
(otherwise,
$c(d,k)$
would be equal to the class of the lines in $D$),
so we have to do some analysis on
the locus where the rank of the homomorphism drops.
The formal proof is given in section 2.

The theorem gives us a tool to compute $\sigma_K$ and $\sigma_L$
using standard techniques in the theory of Chern classes
and Schubert calculus.
Using the symbolic program Maple,
we computed some examples.
The results are very interesting and
a partial list is given in the end of section 2.
In all the cases,
we see that
the union of $\sigma_K$ and $\sigma_L$ is exactly the set of
limiting lines in terms of classes.
This turns out to be true in general.
In fact,
our second main result states that:

\proclaim{Theorem 3.1}
As a class, the sum of $\sigma_K$ and $\sigma_L$
is equal to the class of the lines on a generic hypersurface
that degenerates into the union of $K$ and $ L$.
\endproclaim

The proof of Theorem 3.1 is given in section 3 which
uses mainly elementary combinatorial methods.
Our proof actually gives a more general formula (Theorem 3.13)
involving the Chern classes of
various symmetric powers for any rank two bundle
which should be interesting in its own right.

Finally,
we give in section 4
a proof that
the union of $\sigma_K$ and $\sigma_L$ is indeed the set of
the limiting lines
(Theorem 4.1).
The proof is mainly a formal generalization of Katz's method [K]
in the case of quintic threefolds and one can find the basic idea there.

The author want to thank Katz for some helpful discussions and to
thank the referee for many excellent suggestions,
including an important correction on example (4.10).

\smallpagebreak
\heading 2. The cycle classes of $\sigma_K$ and $\sigma_L$
\endheading
\smallpagebreak

We will consider $\{X_s\}_{s\in \Delta}$
a semistable degeneration of generic hypersurfaces of degree $d$
in $\Bbb P^n$ such
that the central fiber $X_0$ has two irreducible components.
%By choosing the parameter space carefully,
%we may also assume that the lines in $X_s$ are parametrized by $s$ over
%$\Delta^*$.
As indicated in the introduction,
we will use the same letter to denote a hypersurface
and its equation.
Up to the first order,
we may written that
$$X_s=KL+sD,$$
where $K$, $L$, and $D$ are generic hypersurfaces of degree $k$, $l$,
and $d=k+l$,
respectively.
Let $G(2,n+1)$ be the Grassmannian of the lines in $\Bbb P^n$.
We define
$$\sigma_s=\{\alpha \in G(2,n+1)| \  \alpha \subset X_s\},$$
$$\sigma_K=\{\alpha \in G(2,n+1)| \  \alpha \subset K, \alpha \cap L \subset
D\},$$
$$\sigma_L =\{\alpha \in G(2,n+1)| \  \alpha \subset L, \alpha \cap K \subset
D\}.$$
Geometrically speaking,
$\sigma_K$ (respectively, $\sigma_L$) is a subset of $\sigma_0$
which consists of lines $\alpha$ in $K$
(respectively, in $L$) subject to the condition
that $\alpha$ intersects $K\cap L\cap D$ at least $l$ times
(respectively, $k$ times).
This extra condition is imposed by the infinitesimal deformation.
In fact,
one sees from Katz's method in [K]:

\proclaim{Lemma 2.1 (Katz, [K])}
Let $\sigma_{lim} \subset \sigma_0$ be the set of the limiting lines.
By considering
the first order obstructions,
one have that
$$\sigma_{lim} \subset \sigma_K \cup \sigma_L .$$
\endproclaim

Notice that both sides in the expression above have the same expected
dimension.
A natural question then is whether the union of
$\sigma_K$ and $\sigma_L$
is actually equal to $\sigma_{lim}$,
as shown by Katz [K] in the case of quintic threefolds.
Our first goal is to compute the cycle classes of those $\sigma$'s in
the Chow ring $A(G(2,n+1))$ of $G(2,n+1)$.
Of course,
it is well-known that
the class of $\sigma_s$,
$s\not = 0$,
is equal to the top Chern class of the vector bundle $S^dU^*$,
where $U$ is the universal bundle on $G(2,n+1)$.
To describe the classes of $\sigma_K$ and $\sigma_L$,
we need first to recall some definitions.
Please refer to Fulton's book [F] for further details.

\proclaim{Definition 2.2}
Let $V$ and $W$ be two vector bundles.
One defines
$$\eqalign{c&(V-W)=c(V)/c(W)=c(V)s(W)=\cr
&1+(c_1(V)-c_1(W))+(c_2(V)-c_1(V)c_1(W)+c_1^2(W)-c_2(W))+\cdots ,\cr}$$
where $c(V)$ is the Chern class of $V$ and
where $s(W)$ is the Segre class (the inverse Chern class) of $W$.
\endproclaim

\proclaim{Theorem 2.3}
Let $c(\sigma_K)$ be the class of $\sigma_K$.
We then have that
$$c(\sigma_K)=c_{k+1}(S^kU^*)c_l(S^{k+l}U^*-S^kU^*) .\tag 2.4$$
The corresponding representation holds for $\sigma_L$.
\endproclaim

\demo{Proof of Theorem 2.3}
Recall that
$$\eqalign{\sigma_K=&\{\alpha \in G(2,n+1)| \  \alpha \subset K,\alpha \cap L
\subset D\}\cr
=& \{\alpha \in G(2,n+1)| \  \alpha \subset K\}
\cap \{\alpha \in G(2,n+1) | \  \alpha \cap L \subset D\}.\cr}$$
Since $K$ is generic,
$$c(\sigma_K)=c_{top}(S^kU^*)c(D,L),$$
where $c(D,L)$ is the class of
$$\{\alpha \in G(2,n+1) | \  \alpha \cap L \subset D\}.$$
To see $c(D,L)$,
we consider the following exact sequence:
$$ S^kU^* @>\rho_L>> S^{k+l}U^* @>>> B_K @>>> 0 ,\tag 2.5$$
where $\rho_L$ is given by tensoring the section of $S^lU^*$ induced by $L$
and $B_K$ is the quotient.
Notice that $\rho_L$ is not injective as a homomorphism of vector bundles.
To see the degeneracy locus of $\rho_L$,
we observe that,
over $\alpha \in G(2,n+1)$,
$\rho_L$ is given by the restriction
on the degree $k$ part
of the multiplication of $L|_{\alpha}$ in the graded algebra
induced from $\oplus_{i=0}^\infty S^iU^*$.
Therefore,
$\rho_L$ drops rank at $\alpha$ if and only if
$L|_{\alpha}=0$,
that is,
if and only if $\alpha \subset L$.
This implies that the class of the degeneracy locus of $\rho_L$ is equal to
$\lambda c_{l+1}(S^lU^*)$,
where $\lambda $ is the multiplicity.
(Since $\lambda $ is not needed for our proof,
we will compute it later.)
Let $\delta$ be the support of this locus.
On $G^*=G(2,n+1)-\delta$,
we have an exact sequence of vector bundles
$$ 0 @>>> S^kU^* @>\rho_L>> S^{k+l}U^* @>>> B_K @>>> 0 .$$
Therefore,
restricting to $G^*$,
$$B_K=S^{k+l}U^*/\rho_L(S^kU^*)$$
is a vector bundle of rank $l$
and
$$c(B_K)=c(S^{k+l}U^*-S^kU^*).$$
Now $D$ induces a section of $S^{k+l}U^*$B
and hence a section of $B_K$.
At each $\alpha \in G^*$,
the value of the section induced by $D$ in $S^{k+l}U^*$ is given by the
intersection
of $\alpha$ with $D$.
It lies in $\rho_L(S^kU^*)$ if and only if this intersection
is the union of the intersection of $\alpha$ with $L$ and
the intersection of $\alpha$ with some degree $k$ hypersurface.
Therefore,
the class of the zero locus of the section induced by $D$
in $B_K$ is exactly $c(D,L)$.
On the other hand,
it is easy to count that $c(D,L)$ has the expected codimension $l$.
Therefore,
the class of the zero locus of the section induced by
$D$ in $B_K$ is equal to the top
Chern class of $B_K$.
In summary,
we have that
on $G^*$
$$c(D,L)=c_{top}(B_K)=c_l(S^{k+l}U^*-S^kU^*).$$
To see why the formula above still holds on $G(2,n+1)$,
we consider the open inclusion
$$j: G^* \to G(2,n+1).$$
This induces an exact sequence
$$A(\delta) @>i_*>> A(G(2,n+1)) @>j^*>> A(G^*) @>>> 0 ,$$
where
$$i: \delta \to G(2,n+1)$$
is the inclusion of $\delta$ in $G(2,n+1)$.
Notice that the codimension of $\delta$ in $G(2,n+1)$ is $l+1$.
Therefore,
$j^*$ is an isomorphism on $A_m$ for $m \leq l$.
For the convenience of further use,
we summarize this simple observation as the following:

\proclaim{Lemma 2.6}
Let $Y \subset X$ be a closed subvariety of
codimension $n$ in $X$.
For any vector bundles $V$ and $W$ on $X$ and $m < n$,
one has that
$$j_*c_m(j^*V)=c_m(V), \qquad j_*c_m(j^*V-j^*W)=c_m(V-W),$$
where
$$j: X-Y \to X$$
is the inclusion.
\endproclaim

The second identity follows from the fact that $c_m(V-W)$
is defined in terms of the Chern classes of $V$ and $W$
in codimensions less than or equal to $m$.

Putting everything together,
we finally have that
$$c(D,L)=j_*c_{top}(j^*B_k)=j_*c_{l}(j^*S^{k+l}U^*-j^*S^kU^*)
=c_{l}(S^{k+l}U^*-S^kU^*).$$
This gives us the representation for $\sigma_K$.
The same argument gives the representation for $\sigma_L$.
\qed
\enddemo

Theorem 2.3 gives us a tool to compute
$\sigma_K$ and $\sigma_L$
using standard techniques in the theory of Chern classes
and Schubert calculus.
We have made some computations using a home--made program in Maple.
To keep the track of the degrees of the hypersurfaces in our degeneration,
we will use $c(d,k)$ for
$c(\sigma_K)$ and $c(d,l)$ for $c(\sigma_L)$
in the following examples.

\proclaim {Example 1} Lines on a generic cubic surface.
\endproclaim
It is well-known that there are 27 lines on such surface.
When the surface degenerates to a quadric $K$ and a plane $L$,
one easily computes from Theorem 2.3 that
$$c(3,2)=12,\quad c(3,1)=15.$$
It might be interesting to see this directly from the geometry.
We see that $K \cap L \cap D$ consists of six points.
Since for each point in quadric $K$ there are exactly two lines in $K$
passing through it,
$\sigma_K$ contains $12$ lines.
On the other hand,
a line is in
$\sigma_L$ if and only if it passes through two of those six points.
Therefore,
$\sigma_L$ indeed contains $15$ lines.

\proclaim {Example 2} Lines on a generic quintic threefold.
\endproclaim
There are 2,875 lines in a generic quintic threefold.
There are two degenerations
and computations using our representations show:
\roster
\itemitem{} Case 1, $K$ a quartic and $L$ a plane:
$\quad c(5,4)=1600$,
$c(5,1)=1275$.
\itemitem{} Case 2, $K$ a cubic and $L$ a quadric:
$\quad c(5,3)=1575$,
$c(5,2)=1300$.
\endroster
Very interesting computations
of the same results using the geometric method
was first given by Katz in [K].

\proclaim {Example 3} Lines on a generic fourfold of degree $7$.
\endproclaim
One computes that there are $698,005$ lines on such fourfold.
There are three degenerations.
Our computations show that:
\roster
\itemitem{} Case 1, $k=6$ and $l=1$:
$\quad c(7,6)=423,360$,
$c(7,1)=274,645$.
\itemitem{} Case 2, $k=5$ and $l=2$:
$\quad c(7,5)=398,125$,
$c(7,2)=299,880$.
\itemitem{} Case 3, $k=4$ and $l=3$:
$\quad c(7,4)=360,640$,
$c(7,3)=337,365$.
\endroster
Notice that in each of the cases above,
$$c(7,k)+c(7,l)=698,005=c_{top}(S^7U^*),$$
the total number of the lines in a generic
fourfold of degree $7$.

%\proclaim {Example 4} Lines on a generic fivefold of degree $9$.
%\endproclaim
%Our computer shows that
%$$c_{top}(S^9U^*)=305,093,061.$$
%For five degenerations,
%the results are as follows:
%\roster
%\itemitem{} Case 1, $k=8$ and $l=1$:
%$\quad c(9,8)=193,665,024$,
%$c(9,1)=111,428,037$.
%\itemitem{} Case 2, $k=7$ and $l=2$:
%$\quad c(9,7)=179,499,789$,
%$c(9,2)=125,593,272$.
%\itemitem{} Case 3, $k=6$ and $l=3$:
%$\quad c(9,6)=166,453,056$,
%$c(9,3)=138,640,005$.
%\itemitem{} Case 3, $k=5$ and $l=4$:
%$\quad c(9,5)=157,227,525$,
%$c(9,4)=147,865,536$.
%\endroster
%We again have that in all the cases the sum of
%$c(9,k)$ and $c(9,l)$ is equal to
%the total number of the lines in a generic
%fourfold of degree $9$.

\proclaim {Example 4} Lines on a generic quadric in $\Bbb P^n$, $n\geq 3$.
\endproclaim
It is easy to see that the class of the lines on a smooth
quadric is $4\sigma_{2,1}$.
(See also [Chapter 7, GH] for computations using the geometric method.)
When the quadric degenerates into two $(n-1)$-planes $K$ and $L$,
our formula shows that
$$c(2,1)=2\sigma_{2,1}.$$
Once again,
we have that the sum of $c(d,k)$ and $c(d,l)$ is equal to the expected class.
It is easy to check directly that $c(2,1)$ is indeed equal to
the class of $\sigma_K$.

In theory,
one can compute any $c(d,k)$ and $c(d,l)$ in $A(G(2,n+1))$.
In practice,
one needs a computer when $d$ and $n$ are large.
The examples above inspired us to suspect that in general
the sum of $c(d,k)$ and $c(d,l)$ should be equal to the top Chern class
of $S^dU^*$.
It turns out that it is indeed the case and
we will show this fact in the next section.
The upshot of all this is
that the union of $\sigma_K$ and $\sigma_L$ should be
exactly equal to $\sigma_{lim}$ of limiting lines.

\smallpagebreak
\heading 3. An identity in difference classes
\endheading
\smallpagebreak

In this section,
we will show that,
as a class,
the union of $\sigma_K$ and $\sigma_L$
is indeed equal to the class of the lines contained in a generic hypersurface
which degenerates to the union of $K$ and $L$.
More precisely,
we will show:

\proclaim{Theorem 3.1}
Retaining the notation in the last section,
we have the following identity:
$$c_{k+1}(S^kU^*)c_l(S^{d}U^*-S^kU^*)
+c_{l+1}(S^lU^*)c_k(S^{d}U^*-S^lU^*)=c_{d+1}(S^{d}U^*).\tag 3.2$$
In another word,
as a class,
the sum of $\sigma_K$ and $\sigma_L$
is equal to $\sigma_{lim}$.
\endproclaim

To show Theorem 3.1,
first notice that,
by the splitting principle,
we may assume that $U^*$ is the direct sum of two line bundles $A$ and $B$.
Let $a$ be the first Chern class of $A$
and $b$ be the first Chern class of $B$.
Using the standard techniques in computing the Chern classes,
we then have that the total Chern class of $S^kU^*$ is given by
$$ c(S^kU^*)= \prod_{i=0}^{k}[1+(k-i)a+ib],$$
that is,
the $j$-th Chern class of $S^kU^*$ is equal to the degree $j$ part of
the symmetric polynomial in $a$ and $b$ as given in the expression above.
Therefore,
it is enough to show the following polynomial identity:
\proclaim{Theorem 3.3}
Let $c_i(k)$ be defined by
$$ \prod_{i=0}^{k}\{1+[(k-i)a+ib]x\}=\sum_{i=0}^{k+1}c_i(k)x^i,
\quad \text{and $c_i(k)=0$ for $i > k+1$}.$$
If one define inductively
$$\left\{ \eqalign{c_0(k+l,l) & =1,\cr
c_i(k+l,l) & =c_i(k+l)-\sum_{j=0}^{i-1} c_j(k+l,l)c_{i-j}(l),
\quad i > 0 ,\cr}\right.$$
then the following identity holds:
$$c_{k+1}(k) c_l(k+l,k) + c_{l+1}(l) c_k(k+l,l)=c_{k+l+1}(k+l).$$
\endproclaim

\demo{Proof of Theorem 3.3}
The proof is elementary and
we will break it into a few lemmas.

To begin with,
using induction,
we see from the definition of $c_n(k+l,l)$ that
$$ c_n(k+l,l)=\sum_{i=0}^nf_i(a,b;l)c_{n-i}(k+l),$$
where $f_i(a,b;l)$ is a symmetric homogeneous polynomial of degree $i$
in $a$ and $b$ with parameter $l$.
Notice that $f_i$ is independent of $k$ and $n$.

\proclaim{Lemma 3.4}
$$f_i(a,b;l)=\sum_{j=0}^l{(-1)^{i+j}[(l-j)a+jb]^{l+i}
\over j!(l-j)!(a-b)^l}.$$
\endproclaim

\demo{Proof of Lemma 3.4}
{}From the definition,
we have that

$$\eqalign{\sum_{i=0}^{\infty}f_i(a,b;l)x^i
=&{1 \over \prod_{j=0}^{l}\{1+[(l-j)a+jb]x\}}\cr
=&\sum_{j=0}^{l}{h_j(a,b;l) \over 1+[(l-j)a+jb]x}\cr
=&\sum_{j=0}^{l}h_j(a,b;l) \sum_{i=0}^{\infty}(-1)^i[(l-j)a+jb]^ix^i\cr
=&\sum_{i=0}^{\infty} \sum_{j=0}^{l}h_j(a,b;l)
(-1)^i[(l-j)a+jb]^ix^i.\cr}$$
Therefore,
$$f_i(a,b;l)
=\sum_{j=0}^{l}h_j(a,b;l) (-1)^i[(l-j)a+jb]^i.$$
To find $h_j(a,b;l)$,
notice that
$$1=\sum_{j=0}^{l}h_j(a,b;l)
\prod_{i=0,i\not = j}^{l}\{1+[(l-i)a+ib]x\}.$$
For any fixed integer $j$ between $0$ and $l$,
let
$$x=-[(l-j)a+jb]^{-1}.$$
We then get that

$$h_j(a,b;l)={(-1)^j[(l-j)a+jb]^l\over (a-b)^lj!(l-j)!}.$$
This gives us the formula for $f_i(a,b;l)$.
\qed
\enddemo

\proclaim{Lemma 3.5}
$$\sum_{j=0}^l {(-1)^{j+l} \over j!(l-j)!(a-b)^l[(l-j)a+jb]}
={1 \over{\prod_{j=0}^l[(l-j)a+jb]}}
={1 \over c_{l+1}(l)}.$$
\endproclaim

\demo{Proof of Lemma 3.5}
We only need to show that
$$\sum_{j=0}^{l}{(-1)^{j+l}\over j!(l-j)!}
\prod_{i=0,i\not = j}^{l}[(l-i)a+ib]=(a-b)^l.\tag{3.6}$$
Notice that both sides of (3.6) are homogeneous
polynomials of degree $l$ in $a$ and $b$.
Therefore,
we need only to check (3.6) for $l+1$ pairs of $a$ and $b$
with distinct quotients $a/b$.
For $h=0,1,\dots, l$,
let $a=h$ and $b=h-l$.
This gives us that
$$\eqalign{\text{LHS of (3.6)}
=&\sum_{j=0}^{l}{(-1)^{j+l}\over j!(l-j)!}
\prod_{i=0,i\not = j}^{l}[(l-i)h+i(h-l)]\cr
=&{(-1)^{h+l}\over h!(l-h)!}
\prod_{i=0,i\not = h}^{l}(lh-il)\cr
=&l^l\cr
=&(a-b)^l.\quad \qed\cr}
$$
\enddemo

\proclaim{Lemma 3.7}
$$\sum_{j=0}^l {(-1)^j[(l-j)a+jb]^n \over j!(l-j)!}=0,
\quad \text{for $0\leq n \leq l-1$}.$$
\endproclaim

\demo{Proof of Lemma 3.7}
Let $a=h$ and $b=h-l$.
We only need to show that
$$\eqalign{
&\sum_{j=0}^l {(-1)^j[(l-j)h+j(h-l)]^n \over j!(l-j)!}\cr
=&\sum_{j=0}^l {(-1)^jl^n(h-j)^n \over j!(l-j)!}\cr
=&\sum_{j=0}^l {(-1)^jl^n \over j!(l-j)!}\sum_{i=0}^n {n \choose i}
h^{n-i}(-j)^i\cr
=&\sum_{i=0}^n {n \choose i} (-1)^ih^{n-i}l^n
\sum_{j=0}^l {(-1)^jj^i \over j!(l-j)!}\cr
=&0 .}$$
Therefore,
it is enough to show that
$$\sum_{j=0}^l {(-1)^jj^i \over j!(l-j)!}=0, \quad \text{for $0\leq i \leq
l-1$}.$$
This can be easily verified by considering the following:
$$\left\{ \eqalign{f_0(x)&=(1+x)^l=\sum_{j=0}^l {l \choose j} x^j;\cr
f_i(x)&=xf_{i-1}', \quad i > 0 .\cr}\right.$$
In fact,
it is easy to see that,
for $0\leq i \leq l-1$,
$$0=f_i(-1)=
\sum_{j=0}^l {l \choose j} (-1)^j j^i=
l!\sum_{j=0}^l {(-1)^j j^i \over j!(l-j)!}.\quad \qed$$
\enddemo

We continue our proof of Theorem 3.3.
By Lemma 3.4,
$$
\eqalign{&c_k(k+l,l)\cr
=&\sum_{i=0}^kf_i(a,b;l)c_{k-i}(k+l)\cr
=&\sum_{i=0}^k\sum_{j=0}^{l}{(-1)^{i+j}[(l-j)a+jb]^{l+i}\over
j!(l-j)!(a-b)^l}c_{k-i}(k+l)\cr
=&\sum_{j=0}^l\sum_{i=0}^k{(-1)^{l+1+j}[-(l-j)a-jb]^{l+i+1} \over
j!(l-j)!(a-b)^l[(l-j)a+jb] }c_{k-i}(k+l) \cr
=&\sum_{j=0}^l{(-1)^{l+1+j}\over (a-b)^lj!(l-j)![(l-j)a+jb] }
\sum_{i=-l-1}^k[-(l-j)a-jb]^{l+i+1} c_{k-i}(k+l)\cr
&-\sum_{j=0}^l{(-1)^{l+1+j}\over (a-b)^lj!(l-j)![(l-j)a+jb] }
\sum_{i=-l-1}^{-1}[-(l-j)a-jb]^{l+i+1} c_{k-i}(k+l)\cr
=&\sum_{j=0}^l{(-1)^{l+1+j}\over (a-b)^lj!(l-j)![(l-j)a+jb]}
\sum_{i=0}^{k+l+1}[-(l-j)a-jb]^i c_{k+l+1-i}(k+l)\cr
&-\sum_{j=0}^l{(-1)^{l+1+j}\over (a-b)^lj!(l-j)![(l-j)a+jb]}
\sum_{i=0}^l[-(l-j)a-jb]^i c_{k+l+1-i}(k+l)\cr
=&\sum_{j=0}^l{(-1)^{l+1+j}
\over j!(l-j)![(l-j)a+jb](a-b)^l}
\prod_{i=0}^{k+l}[-(l-j)a-jb+(k+l-i)a+ib]\cr
&+\sum_{i=0}^l (-1)^ic_{k+l+1-i}(k+l)
\sum_{j=0}^l{(-1)^{l+j}[(l-j)a+jb]^{i-1}\over (a-b)^lj!(l-j)!}.\cr}$$

Therefore,
using Lemma 3.5 and Lemma 3.7,
we have the following:
$$c_k(k+l,l)=-c_{k+1}(k)f+c_{k+l+1}(k+l)/c_{l+1}(l),\tag 3.8$$
where
$$f=\sum_{j=0}^l{(-1)^{l+j}
\prod_{i=0}^{j-1}[(k+j-i)a+(i-j)b]\prod_{i=k+1+j}^{k+l}[(k+j-i)a+(i-j)b]
\over j!(l-j)![(l-j)a+jb](a-b)^l}.$$

To finish our proof of Theorem 3.3,
we need to show the following:

\proclaim{Lemma 3.9}
$$c_{l+1}(l)f=c_l(k+l,k). \tag3.10$$
\endproclaim

\demo{Proof of Lemma 3.9}
To start with,
notice that $c_{l+1}(l)f$ is a homogeneous polynomial of degree $l$
in $a$ and $b$.
In fact,
from (3.8),
$$ c_{k+1}(k)c_{l+1}(l)f=c_{k+l+1}(k+l)-c_k(k+l,l)c_{l+1}(l)$$
is a homogeneous polynomial in $a$ and $b$.
Moreover,
$(a-b)$ is not a factor of $c_{k+1}(k)$.
Therefore,
$c_{l+1}(l)f$ must be a polynomial in $a$ and $b$.
We hence need only to check (3.10) for $a=h$ and $b=l-h$,
$h=0, 1, \dots, l$.
For such $a$ and $b$,
one easily computes that
$$c_{l+1}(l)f={\prod_{i=0}^{k+l}(kh+lh-il)\over
\prod_{i=0}^{k}(kh-il)},
\quad (i\not = kh/l \quad \text{if $kh/l$ is an integer).}$$
On the other hand,
using (3.8) and interchanging $l$ and $k$,
we see that
$$c_l(k+l,k)=-c_{l+1}(l)g+c_{k+l+1}(k+l)/c_{k+1}(k),$$
where
$$g=\sum_{j=0}^k{(-1)^{k+j}
\prod_{i=0}^{j-1}[(l+j-i)a+(i-j)b]\prod_{i=l+1+j}^{k+l}[(l+j-i)a+(i-j)b]
\over j!(k-j)![(k-j)a+jb](a-b)^k}.$$
For $a=h$, $b=h-l$, and $h=0,1.\dots, l$,
we see that $c_{l+1}(l)=0$.
Therefore,
if
$$(k-j)h+j(h-l)\not = 0,$$
that is,
if $kh/l$ is not an integer,
we have that
$$c_{l+1}(l)g=0.$$
This gives us that
$$c_l(k+l,k)=c_{k+l+1}(k+l)/c_{k+1}(k)
={\prod_{i=0}^{k+l}(kh+lh-il)\over
\prod_{i=0}^{k}(kh-il)}=c_{l+1}(l)f.$$
If $kh/l=n$ is an integer,
then we take
$$a=h+\delta, \quad b=h+\delta-l.$$
It is a straightforward computation that
$$\eqalign{\lim_{\delta\to 0} c_l(k+l,k)
=&-\lim_{\delta\to 0}c_{k+1}(k)g+\lim_{\delta\to 0}c_{k+l+1}(k+l)/c_{k+1}(k)\cr
=&-{l\over k}{\prod_{i=0,i\not = n}^{k+l}(kh+lh-il)
\over\prod_{i=0,i\not = n}^{k}(kh-il)}
+{k+l\over k}{\prod_{i=0,i\not =
n}^{k+l}(kh+lh-il)\over\prod_{i=0,i\not = n}^{k}(kh-il)}\cr
=&c_{l+1}(l)f .\cr}$$
This completes our proof of Lemma 3.9, Theorem 3.3, and hence Theorem 3.1.
\qed
\enddemo
\enddemo

As we mentioned in section 2,
the degeneracy locus of $\rho_L$ in (2.5) is equal to
$\lambda c_{l+1}(S^lU^*)$,
where $\lambda $ is the multiplicity.
We will now show this fact.

\proclaim{Proposition 3.11}
Let $\rho_L$ be the homomorphism in 2.5 and
$\Bbb D(\rho_L)$ be the cycle class of the degeneracy locus
$$D(\rho_L)=\{\alpha \in G(2,n+1)| \  \text{rank}(\rho_L(\alpha))\leq k\},$$
where $\rho_L(\alpha)$ is the restriction of $\rho_L$ on the fiber over
$\alpha$.
We then have
$$\Bbb D(\rho_L)=\lambda c_{l+1}(S^lU^*),$$
where $\lambda = \sum_{j=0}^k{l+j \choose j}$.
\endproclaim

\demo{Proof of Proposition 3.11}
By the analysis on the degeneracy locus of $\rho_L$
in the proof of Theorem 2.3,
we have
$$\Bbb D(\rho_L)=\lambda c_{l+1}(S^lU^*),$$
and need only to find the multiplicity $\lambda$.
On the other hand,
since $D(\rho_L)$ has the expected codimension,
by the Thom--Porteous formula [Chapter 14, F],
we have that
$$\Bbb D(\rho_L)=c_{l+1}(S^{k+l}U^*-S^kU^*).$$
This gives an easy way to compute $\lambda$ and that will be the end of
the proof.
However,
it might be more convincing to show the following lemma directly.

\proclaim{Lemma 3.12}
$$c_{l+1}(S^{k+l}U^*-S^kU^*)=\lambda c_{l+1}(S^lU^*),$$
where $\lambda = \sum_{j=0}^k{l+j \choose j}$.
\endproclaim

\demo{Proof of Lemma 3.12}
Using the notation in the proof of Theorem 3.1,
we need only to show the polynomial identity
$$c_{l+1}(k+l,k)=\lambda c_{l+1}(l).$$
By Lemma 3.4 (replacing $l$ by $k$),
$$\eqalign{c_{l+1}(k+l,k)
=&\sum_{i=0}^{l+1}f_i(a,b;k)c_{l+1-i}(k+l)\cr
=&\sum_{i=0}^{l+1}\sum_{j=0}^{k}{(-1)^{i+j}[(k-j)a+jb]^{k+i}\over
j!(k-j)!(a-b)^k}c_{l+1-i}(k+l)\cr
=&\sum_{j=0}^k\sum_{i=0}^{l+1}{(-1)^{k+j}[-(k-j)a-jb]^{k+i} \over
j!(k-j)!(a-b)^k}c_{l+1-i}(k+l) \cr
=&\sum_{j=0}^k{(-1)^{k+j}\over (a-b)^kj!(k-j)!}
\sum_{i=-k}^{l+1}[-(k-j)a-jb]^{k+i} c_{l+1-i}(k+l)\cr
&-\sum_{j=0}^k{(-1)^{k+j}\over (a-b)^kj!(k-j)!}
respectively,
\sum_{i=-k}^{-1}[-(k-j)a-jb]^{k+i} c_{l+1-i}(k+l)\cr
=&\sum_{j=0}^k{(-1)^{k+j}\over (a-b)^kj!(k-j)!}
\sum_{i=0}^{k+l+1}[-(k-j)a-jb]^i c_{k+l+1-i}(k+l)\cr
&-\sum_{i=0}^{k-1}c_{k+l+1-i}(k+l)
\sum_{j=0}^k{(-1)^{k+j}[-(k-j)a-jb]^i\over (a-b)^kj!(k-j)!}\cr
=&\sum_{j=0}^k{(-1)^{k+j}\over (a-b)^kj!(k-j)!}
\sum_{i=0}^{k+l+1}[-(k-j)a-jb]^i c_{k+l+1-i}(k+l).\cr}$$

In the last step,
we have used Lemma 3.7 (replacing $l$ by $k$)that
the second term in the next to the last step of the expression above is zero.
Using this expression of $c_{l+1}(k+l,k)$,
it is easy to see that,
for $a=h$ and $b=h-l$, $h=0,1, \dots, l$,
$$c_{l+1}(k+l,k)=0.$$
The same is true for $c_{l+1}(l)$.
Since both $c_{l+1}(k+l,k)$ and $c_{l+1}(l)$ are homogeneous
polynomial of degree $l+1$ in $a$ and $b$,
the quotient
$c_{l+1}(k+l,k)/ c_{l+1}(l)$ must be a constant.
To find this constant $\lambda$,
we let $a=1$ and $b \to 0$.
Using our expression of $c_{l+1}(k+l,k)$ again,
one easily computes that $\lambda = \sum_{j=0}^k{l+j \choose j}$.
\qed
\enddemo
\enddemo

We have actually proved more.
Since the proofs in this section did not use any special properties of
the Grassmannian and its universal bundle,
our results can be generalized
to the following forms for any rank two bundle.

\proclaim{Theorem 3.13}
For a vector bundle $V$ of rank two and any pair of positive integers $k$ and
$l$,
the following formula holds:
$$c_{k+1}(S^kV)c_l(S^{k+l}V-S^kV)+c_{l+1}(S^lV)c_k(S^{k+l}V-S^lV)
=c_{k+l+1}(S^{k+l}V).$$
\endproclaim

\proclaim{Proposition 3.14}
The following identity holds
for a vector bundle $V$ of rank two and any pair of positive integers $k$
and $l$:
$$c_{l+1}(S^{k+l}V-S^kV)=\lambda c_{l+1}(S^lV),$$
where $\lambda = \sum_{j=0}^k{l+j \choose j}$.
\endproclaim

\demo{Question}
Is it possible for Theorem 3.13 to be generalized
to vector bundles of higher ranks?
If such formula does exist,
what would it imply geometrically?
\enddemo

The question should have some bearing on degenerations of hypersurfaces
and their $(m-1)$-linear spaces,
where $m$ is the rank of the bundle in consideration.
It is possible,
of course,
that finding the would-be formula might be
as difficult
as actually proving the formula.
For this reason,
a geometric point of view could be very important.

\smallpagebreak
\heading 4. The union of $\sigma_K$ and $\sigma_L$ is $\sigma_{lim}$
\endheading
\smallpagebreak

Theorem 2.3 and Theorem 3.1 strongly suggest that the union of
$\sigma_K$ and $\sigma_L$ should in fact be equal to $\sigma_{lim}$.
We will now prove this fact following Katz's method in [K].
One can find the basic idea and omitted details there.
We will use the same notation as given in section 2.
Although the results obtained so far are valid for all positive integers
$l$, $k$, and $d$ with $k+l=d$,
we will in this section require that those integers to be less than or
equal to $2n-3$,
since a generic hypersurface in $\Bbb P^n$ will possess lines only if
its degree is less than or equal to $2n-3$.

\proclaim{Theorem 4.1} $$\sigma_{lim} = \sigma_K \cup \sigma_L .$$
\endproclaim

We will prove Theorem 4.1 by formally generalizing Katz's method
as given in the case of quintic threefold.
For a rational curve $C$ of degree
$m$ in $\Bbb P^n$,
one can identify $C$ with
its parameter equation
as an element in $\oplus^{n+1}H^0(\Cal O_{\Bbb P^1}(m))$.
Using a slightly different notation from that used in [K],
for any given line $\alpha=(\alpha_0,\alpha_1,\hdots,\alpha_n)$ in $\Bbb P^n$,
one defines a linear map
$$\Phi_{K,\alpha} : \oplus^{n+1}H^0(\Cal O_{\alpha}(1))
\to H^0(\Cal O_{\alpha}(k))$$
by using homogeneous coordinates to set
$$ \Phi_{K,\alpha}(a_0,a_1,\hdots,a_n)=
\sum_{i=0}^na_i\frac{\partial K}
{\partial X_i}(\alpha_0,\alpha_1,\hdots,\alpha_n).$$
Furthermore,
let $\Psi_{D,K,L,\alpha}$
be the composition
$$\Psi_{D,K,L,\alpha} :
\text
{Ker$(\Phi_{K,\alpha}) \subset \oplus^{n+1}H^0(\Cal O_{\alpha}(1))
@>\Phi_{D,\alpha}-\frac{D}{L}\Phi_{L,\alpha}>>
H^0(\Cal O_{\alpha}(d))$}
\to H^0(\Cal O_{\alpha}(d))|_{\alpha\cap L},$$
that is,
the domain of
$\Psi_{D,K,L,\alpha}$ is $\text{Ker}(\Phi_{K,\alpha})$ and the target
of the map is restricted to $l$ points
$\alpha\cap L$.

\proclaim{Theorem 4.2 (Katz)}
In a generic degeneration as given in section 2,
for any $\alpha$ in $\sigma_K$,
$\alpha$ is contained in $ \sigma_{lim}$ if following two conditions are
satisfied:
\endproclaim

\proclaim{Condition 4.3}
Map $\Phi_{K,\alpha}$ is surjective.
\endproclaim

\proclaim{Condition 4.4}
Map $\Psi_{D,K,L,\alpha}$ is surjective.
\endproclaim

For a proof of Theorem 4.2,
please see [K].

We will now verify that two conditions in Theorem 4.2 hold
for a generic $\alpha$ in $\sigma_K$.
This will be enough to show Theorem 4.1,
since we have already proved that
the union of $\sigma_K$ and $\sigma_L$
is equal to $\sigma_{lim}$
as classes.

We will first prove the following stronger version of condition 4.3:

\proclaim{Lemma 4.5}
For a generic hypersurface $K$,
$\Phi_{K,\alpha}$ is surjective for every $\alpha$ in $K$.
\endproclaim

\demo{Proof of Lemma 4.5}
We will consider two commutative diagrams with exact rows
and columns as Clemens [C]
did in the case of quintic threefolds:

$$\CD
@. 0 @. 0 @. @.\\
@. @VVV @VVV @. @.\\
@. \Cal O_{\alpha} @. \Cal O_{\alpha} @. @.\\
@. @VVV @VVV @.\\
0 @>>> \tilde\Cal T_K \otimes \Cal O_{\alpha}
@>>> \Cal O_{\alpha}^{n+1}(1) @>\Phi_{K,\alpha}>> \Cal O_{\alpha}(k)
@>>> 0 \\
@. @VVV @VVV @. @.\\
0 @>>> \Cal T_K \otimes \Cal O_{\alpha}
@>>> \Cal T_{\Bbb P^{n+1}}\otimes \Cal O_{\alpha} @>>>
\Cal O_{\alpha}(k) @>>> 0\\
@. @VVV @VVV @. @.\\
@. 0 @. 0 @. @.
\endCD \tag 4.6
$$

and

$$\CD
@. 0 @. 0 @. @.\\
@. @VVV @VVV @.@.\\
@. \Cal O_{\alpha} @. \Cal O_{\alpha} @. @.\\
@. @VVV @VVV @. @.\\
0 @>>> \tilde\Cal T_{\alpha} @>>> \tilde\Cal T_K\otimes \Cal O_{\alpha}
@>>> \Cal N_{\alpha/K} @>>> 0\\
@. @VVV @VVV @. @.\\
0 @>>> \Cal T_{\alpha} @>>> \Cal T_K\otimes \Cal O_{\alpha}
@>>> \Cal N_{\alpha/K} @>>> 0\\
@. @VVV @VVV @. @.\\
@. 0 @. 0 @. @.
\endCD\tag 4.7
$$

Since $H^1(\Cal O_{\alpha}(1))=0$,
one sees from (4.6) that
$\Phi_{K,\alpha}$ is surjective if and only if
$$H^1(\tilde\Cal T_K\otimes\Cal O_{\alpha})=0.\tag 4.8$$
On the other hand,
since
$$H^1(\tilde\Cal T_{\alpha})
=H^1(\Cal T_{\alpha})=0,$$
one sees from (4.7)
that (4.8) hold if and only if
$$H^1(\Cal N_{\alpha/K})=0. \tag 4.9$$
It is easy to see from the standard exact sequence

$$ 0 \to \Cal N_{\alpha/K}
\to \Cal N_{\alpha/\Bbb P^n}
\to \Cal N_{K/\Bbb P^n} \otimes \Cal O_{\alpha} \to 0
$$
that
$$h^0(\Cal N_{\alpha/K})-h^1(\Cal N_{\alpha/K})=2n-3-k.$$
Therefore,
we need only to check that for a generic $K$
it contains only such $\alpha$'s with
$$h^0(\Cal N_{\alpha/K})=2n-3-k.$$

To see this,
let us consider
$\Cal M=\{ (\alpha, K)| \ \alpha \subset K\} \subset G(2, n+1)\times
\Bbb P^N$.
The fiber over any $\alpha$ in $G(2,n+1)$ is
$\Bbb P(H^0(\Cal I_{\alpha/\Bbb P^n}(k)))$.
{}From the standard exact sequence

$$ 0 \to \Cal I_{\alpha/\Bbb P^n}(k)
\to \Cal O_{\Bbb P^n}(k)
\to \Cal O_{\alpha}(k) \to 0,
$$
it is an easy computation that the fiber has a constant codimension
of $k+1$ in $\Bbb P^N$.
Therefore,
$\Cal M$ is smooth and the fiber over a generic
$K$ in $\Bbb P^N$ has dimension
$$2(n-1)-(k+1)=2n-3-k.$$
Moreover,
by Sard's theorem,
there exists $K$ in $\Bbb P^N$ such that the fiber over $K$ is smooth.
This means that
for a generic hypersurface $K$ of $\Bbb P^n$
the Hilbert scheme $\Cal M_K$ of
$\alpha$'s contained in $K$ is
smooth with the expected dimension.
Since $H^0(\Cal N_{\alpha/K})$
is the tangent space to $\Cal M_K$ at $\alpha$,
this completes our proof.
\qed
\enddemo

\demo{Remarks}

1). The smoothness of $\Cal M_K$ for a generic $K$
also follows from
the fact that $\Cal M_K$ is just
the zero locus of the section of $S^kU^*$
incduced by $K$.

2). The proof of Lemma 4.5 also gives a restriction on possible
types of $\Cal N_{\alpha/K}$ for a generic $K$.
Let $\Cal N_{\alpha/K}= \oplus_{i=1}^{n-2}\Cal O_{\alpha}(a_i)$.
(4.9) holds if and only if all $a_i \geq -1$.
Using this and the facts that $a_i \leq 1$ and $\sum a_i=n-1-k$,
one can determine possible $\Cal N_{\alpha/K}$'s for a generic $K$.
For example,
in the case of $k=2n-3$ and the case of $k=2n-4$,
one sees that a generic $K$
contains only lines with the normal bundle of
type $(-1, -1, \hdots, -1)$
and type $(0, -1, \hdots, -1)$,
respectively.
\enddemo

To verify condition 4.4,
let $\Cal G=\{(D,K,L,\alpha)| \  \alpha \in \sigma_K\}$.
By considering a series of fibrations as Katz did in [K],
one can see that $\Cal G$ is irreducible.
Therefore,
it is enough for us to find one particular point of $\Cal G$
such that condition 4.4 is satisfied and
$\Phi_{K,\alpha}$ is surjective.
The following example is suggested by the referee:

$$\eqalign{
D : \quad &x_nx_0^{d-1}+x_{n-1}x_0^{d-3}x_1^{2}
+\hdots+x_{n-r+1}x_0^{d-2r+1}x_1^{2r-2}=0,\cr
K : \quad &x_2x_0^{k-1}+x_3x_0^{k-3}x_1^{2}
+\hdots+x_{m+1}x_0^{k-2m+1}x_1^{2m-2}+x_{m+2}x_1^{k-1}
=0,\cr
L : \quad &x_0^l-x_1^l=0,\cr
\alpha : \quad &x_2=x_3=\hdots=x_n=0,\cr}
\tag{4.10}$$
where $m=[k/2]$ and $r=[(l+1)/2]$.
Notice that $d \leq 2n-3$ and
we may assume $k < d$.
Therefore,
$$\text{$m+2 \leq n-1$, $\ $ if $k$ is odd,}$$
$$\text{$m+2 \leq \ \ n\ $,$\ \ $ if $k$ is even.}$$
It is easy to check that $\Phi_{K,\alpha}$ is surjective and
$\text{Ker}(\Phi_{K,\alpha})$ has the minimal
dimension $2n-k+1$.
In fact,
$$\text{Ker}(\Phi_{K,\alpha})=\{(a_0,a_1,0,\hdots,0,a_{m+1},
a_{m+2},a_{m+3},\hdots,a_n)\},$$
where $a_i$'s are linear forms in $x_0$ and $x_1$ subject to conditions
$$\text{$\ \ a_{m+1}=0\ \ $, $\ \ a_{m+2}=0\ \ \ $,
$\ \ $ if $k$ is odd,}$$
$$\text{$a_{m+1}=\mu x_1$, $a_{m+2}=-\mu x_0$,
$\ \ $ if $k$ is even.}$$
On the other hand,
$$\Psi_{D,K,L,\alpha}(a_0,\hdots,a_n)
=\Phi_{D,\alpha}(a_0,\hdots,a_n)|_{\alpha\cap L}.$$
Moreover,
it is easy to see that
$n-r+1 \geq m+2$ with the equality holds if and only if $l+k=2n-3$,
$l$ is odd, and $k$ is even.
Therefore,
$$\{(\mu_1x_0^{d}+\mu_2x_0^{d-1}x_1+
\hdots+\mu_lx_0^{d-l+1}x_1^{l-1})|_{\alpha\cap L}\}
\subseteq \text{Im}(\Psi_{D,K,L,\alpha}).$$
Notice that $l$ points $\alpha \cap L$ are
$\{[1,\xi_i]\}_{1\leq i\leq l}$,
where $\xi_i$'s are $l^{th}$ roots of the unity.
{}From this,
we have
$$\text{Im}(\Psi_{D,K,L,\alpha})=
\text{Im}(A),$$
where $A$ is a $l$ by $l$ matrix with
$a_{i,j}=\xi_i^{j-1}$.
Since
$$\det(A)=\prod_{i>j}(\xi_i-\xi_j)\not = 0,$$
$\Psi_{D,K,L,\alpha}$ is indeed surjective.
This completes our proof of Theorem 4.1.

In the case of the top degree,
both $\sigma_K$ and $\sigma_L$ are finite sets.

\proclaim{Corollary 4.11}
In the case of $d=2n-3$,
there is no monodromy on the family of lines in a generic degeneration
as given in section 2.
\endproclaim

\demo{Remark}
Varley has pointed out to the author that there is a way to
prove Corollary 4.11 directly.
This then implies Theorem 4.1 in the top degree case of $d=2n-3$.

\enddemo

\vfill
\eject

\enddocument
\bye